\documentclass[b5paper,10pt,abstractoff,DIV=calc,headings=normal,twoside]{scrartcl}
\usepackage[english]{babel}
\usepackage{lmodern}
\usepackage{microtype}
\usepackage{amsmath,amssymb,amsthm,amsfonts,graphicx,etoolbox,scrlayer-scrpage}
\usepackage[noblocks]{authblk}
\makeatletter
\patchcmd{\@maketitle}{\huge}{\Large}{}{}
\patchcmd{\abstract}{\quotation}{}{}{}
\AtBeginEnvironment{abstract}{\noindent\ignorespaces}
\AtEndEnvironment{abstract}{\par\mbox{}}
\newcommand{\shortauthor}{}
\newcommand{\shorttitle}{\@title}
\makeatother
\setkomafont{author}{\small}
\setkomafont{title}{\rmfamily\bfseries}
\setkomafont{disposition}{\rmfamily\bfseries}

\def\AMS#1{\par\noindent \textbf{AMS subject classification: }#1\par}
\newcommand{\acknowledgements}{\par\mbox{}\par\noindent\textbf{Acknowledgements: }}
\newcommand{\keywords}[1]{\par\noindent\textbf{Keywords: }#1}
\ofoot[]{}\ifoot[]{}\recalctypearea
\theoremstyle{plain}
\newtheorem{theorem}{Theorem}

\theoremstyle{definition}
\newtheorem{definition}{Definition}

\theoremstyle{remark}

\usepackage{comment}
\usepackage{url}
\usepackage{mathrsfs}
\usepackage{euscript}
\usepackage{mathtools}
\usepackage{amssymb}
\usepackage{amsfonts}
\usepackage[shortlabels]{enumitem}
\usepackage{xcolor}
\usepackage{bm}
\usepackage{subfig}
\usepackage{stmaryrd}

\newcommand{\rangleL}{{L^2(\mathbb{M}^d)}}
\newcommand{\rangleLL}{{L^2(\mathbb{M}^d \times \mathbb{M}^d)}}

\newcommand{\rangleHH}{{\mathbb{H}_p}} %\mathcal{H}_p \otimes \mathcal{H}_p
\newcommand{\Hp}{{\mathcal{H}_p}}
\newcommand{\Hq}{{\mathcal{H}_q}}
\newcommand{\Hprod}{{\mathbb{H}_p}}
\newcommand{\pen}{{\eta}}

\newcommand{\W}{W}

\newcommand{\pit}{{\Pi_2(p,q)}}

\newcommand{\K}{\Lambda}
\newcommand{\dime}{\operatorname{dim}(\mathcal{Y}_\ell)}
\newcommand{\dimep}{\operatorname{dim}(\mathcal{Y}_{\ell'})}

\renewenvironment{abstract}{\bigskip\noindent\begin{minipage}{\textwidth}\setlength{\parindent}{15pt}\paragraph{Abstract:}}{\end{minipage}}

\begin{document}

%%%%%%%%%%%%%%%%%%%%%%%%%%%%%%%%%%%%%%%%%%%%%%%%%%%%%%%%%%%%%%%%%%%%%%%%%%%%%%%%
%%%%%%%%% Please make sure to leave space for page numbers             %%%%%%%%%
%%%%%%%%% There is a placeholder PPP generated in the document to help %%%%%%%%%
%%%%%%%%%%%%%%%%%%%%%%%%%%%%%%%%%%%%%%%%%%%%%%%%%%%%%%%%%%%%%%%%%%%%%%%%%%%%%%%%

\renewcommand{\shortauthor}{A. Caponera}

\title{On the Estimation of Anisotropic Covariance Functions on Compact Two-Point Homogeneous Spaces}

\author[1]{Alessia Caponera\thanks{acaponera@luiss.it}}
\affil[1]{LUISS Guido Carli, Rome, Italy}

\date{}

\maketitle

\begin{abstract}
In this paper, the asymptotic theory presented in \cite{ejs} for spline-type anysotropic covariance estimator on the 2-dimensional sphere is generalized to the case of connected and compact two-point homogeneous spaces.
\end{abstract}


\smallskip

 \AMS{Primary 62G08; secondary 62M.} % please visit http://www.ams.org/mathscinet/msc/msc2010.html

%%%%%%%%%%%%%%%%%%%%%%%%%%%%%%%%%%%%%%%%%%%%%%%%%%%%%%%%%%%%%%%%%%%%%%%%%%%%%%%%

%\begin{figure}
%\includegraphics[width=\textwidth]{samplefig01}
%\caption{Figures should be included in PDF}\label{samplefig01}
%\end{figure}

\section{Introduction}

Let $\mathbb{M}^d$ denote a connected and compact two-point homogeneous Riemannian manifold of dimension $d$. According to \cite{wang}, $\mathbb{M}^d$ belongs essentially to one of the following categories: the unit (hyper)sphere $\mathbb S^d$, $d=1,2,3,\dots$, the real projective spaces $\mathbb P^d(\mathbb R)$, $d=2,3,\dots$, the complex projective spaces $\mathbb P^d(\mathbb C)$, $d=4,6,\dots$, the quaternionic projective spaces $\mathbb P^d(\mathbb H)$, $d = 8, 12, \dots$, and the Cayley projective plane $\mathbb P^d(\operatorname{Cay})$, $d = 16$.

In this paper, we extend the asymptotic theory presented in \cite{ejs} for spline-type anysotropic covariance estimator on the 2-dimensional sphere to the case of connected and compact two-point homogeneous spaces. Specifically, we have a zero-mean second-order random field $X=\{X(u), \, u \in \mathbb{M}^d\}$ that is a random element of $L^2(\mathbb{M}^d)$ with covariance function
$$
C(u,v)= \mathbb{E}[X(u) X(v)],
$$
for $u,v \in \mathbb{M}^d$. We consider also $X_1,\dots,X_n$ independent replicates of $X$ and for the $i$-th replicate we make noisy measurements at $r_i$ random locations on $\mathbb{M}^d$, so that 
$$
\W_{ij} = X_i(U_{ij}) + \epsilon_{ij}, \qquad  i =1,\dots,n, \, j =1,\dots,r_i,
$$
where the $U_{ij}$'s 
%are independent measurement locations uniformly distributed over the sphere $\mathbb{M}^d$
are independently drawn from a common distribution on $\mathbb{M}^d$, and the $\epsilon_{ij}$'s are independent and identically distributed measurement errors of mean 0 and variance $0<\sigma^2<\infty$. %Furthermore, the $X_i$'s and the measurement locations are assumed to be independent of the measurement errors.
Furthermore, the $X_i$'s, the measurement locations, and the measurement errors are assumed to be mutually independent. We address the functional data analysis problem of estimating their second-order moment structure.

In Section 2, we provide an introduction to several key concepts, including the Laplace-Beltrami operator, addition formula, and Sobolev spaces, defined within the context of compact two-point homogeneous spaces. These tools are then used to carry out the asymptotic theory for the spline-type covariance estimator presented in Section 3.

\section{Mathematical Background}\label{sec:background}

Here, the space $\mathbb{M}^d$ is endowed with the usual Riemannian (geodesic) distance $\rho$, and we assume that $\rho(u,v) \in [0,\pi]$, for every $u,v \in \mathbb{M}^d$. We also consider a \emph{normalized} Riemannian measure, denoted by $d\nu$, so that $\int_{\mathbb{M}^d} d\nu = 1$.

Let $\Delta_{\mathbb{M}^d}$ be the Laplace–Beltrami operator on $\mathbb{M}^d$. It is well known that the spectrum of $\Delta_{\mathbb{M}^d}$ is purely discrete, the eigenvalues being
     \begin{equation*}
       \lambda_\ell =  \lambda_\ell (\mathbb{M}^d)  = -\ell (\ell +\alpha +\beta +1),\quad \ell \in \K,
     \end{equation*}
    where
    $$
    \K =  \K (\mathbb{M}^d) = \begin{cases}
   \{n \in \mathbb{N}_0: n \text{ even}\} &\text{if }  \mathbb{M}^d = \mathbb{P}^d(\mathbb R)\\
   \mathbb{N}_0 &\text{if }  \mathbb{M}^d \ne \mathbb{P}^d(\mathbb R),
    \end{cases}
    $$
    and the parameters $\alpha, \beta$ are reported in Table \ref{tab::alphabeta}.
    Now, consider the space of real-valued square-integrable functions over $\mathbb{M}^d$, written $L^2(\mathbb{M}^d)$, endowed with the standard inner product
$$
\langle f, g \rangle_{\rangleL} = \int_{\mathbb{M}^d} f(u)g(u) d\nu(u), \qquad f,g \in L^2(\mathbb{M}^d).
$$

\begin{table}[ht]
\caption{}
    \begin{center}
    \def\arraystretch{1.5}
    \begin{tabular}{c|ccccc}\hline
     $\mathbb{M}^d$ & $\mathbb S^d$ & $\mathbb{P}^d(\mathbb R)$ & $\mathbb{P}^d(\mathbb C)$ & $\mathbb{P}^d(\mathbb H)$ & $\mathbb P^d(\operatorname{Cay})$ \\ \hline
   Dimension& $d=1,2,\dots$ & $d=2,3,\dots$ & $d=4,6,\dots$ & $d=8,12,\dots$ & $d=16$ \\
   $\alpha$ & $(d-2)/2$ & $(d-2)/2$ & $(d-2)/2$ & $(d-2)/2$ & 7 \\
   $\beta$ & $(d-2)/2$ & $(d-2)/2$ &0& 1 & 3 \\
  \hline    \end{tabular}\label{tab::alphabeta}
    \end{center}
\end{table}

    For every $\ell \in \Lambda$, the eigenfunctions corresponding to the same eigenvalue $\lambda_\ell$ form a finite dimensional vector space, denoted by $\mathcal{Y}_\ell = \mathcal{Y}_\ell (\mathbb{M}^d)$, with dimension $\dime$ that is uniquely identified through
    \begin{equation*}
\dime = \frac{(2\ell + \alpha + \beta + 1) \Gamma(\beta +1) \Gamma(\ell +\alpha + \beta +1)\Gamma(\ell +\alpha +1)}{\Gamma(\alpha +1) \Gamma(\alpha + \beta +2)\Gamma(\ell +1)\Gamma(\ell +\beta+1)}. 
\end{equation*}
   Then, 
\begin{equation*}
    L^2(\mathbb M^d) = \bigoplus_{\ell \in \Lambda} \mathcal{Y}_\ell.
\end{equation*}
Moreover, given an orthonormal basis of $\mathcal{Y}_\ell$, say $\{ Y_{\ell,m}, m=1,\dots, \dime \}$, any function $f \in L^2(\mathbb{M}^d)$ can be expanded as 
\begin{equation*}
    f = \sum_{\ell \in \Lambda} \sum_{m =1}^{\dime} \langle f , Y_{\ell,m} \rangle_\rangleL Y_{\ell,m}.
\end{equation*}
In order to deal with covariance functions, we will also consider functions $g \in L^2(\mathbb{M}^d \times \mathbb{M}^d)$, for which it holds
\begin{equation*}
   g= \sum_{\ell, \ell' \in \Lambda} \sum_{m =1}^{\dime}  \sum_{m' =1}^{\dimep}  \langle g , Y_{\ell,m} \otimes Y_{\ell', m'} \rangle_\rangleLL Y_{\ell,m} \otimes Y_{\ell',m'}.
\end{equation*}

A crucial result of \cite[Theorem 3.2]{gine} ensures that also for general $\mathbb{M}^d$ an \emph{addition formula} holds, that is, for any $\ell \in \Lambda,$
\begin{equation*}
\sum_{m=1}^{\dime } Y_{\ell, m}(u) Y_{\ell,m}(v) = \kappa_{\ell}  P_{\ell }^{(\alpha, \beta)}(\cos \epsilon \rho(u,v)),
\end{equation*}
where $\epsilon=1/2$ if $\mathbb{M}^d = \mathbb{P}^d(\mathbb R)$, $\epsilon=1$ otherwise, 
\begin{equation*}
    \kappa_\ell = \kappa_{\ell}({\mathbb{M}^d}) %= \frac{\Gamma(\beta+1)(2\ell +\alpha+\beta+1)\Gamma(\ell +\alpha+\beta+1)}{\Gamma(\alpha +\beta+2)\Gamma(\ell +\beta+1)} 
    = \frac{\dime}{P_{\ell}^{(\alpha, \beta)}(1)}. 
    \end{equation*}      
   The functions $P_{\ell }^{(\alpha, \beta)}:[-1,1]\to \mathbb{R}, \ \ell \in \mathbb{N}_0,$ are the well-known \emph{Jacobi polynomials}. Note that when $\mathbb{M}^d = \mathbb{S}^2$ they reduce to the \emph{Legendre polynomials}.

In general, the reason for the exceptional behaviour of the real projective spaces is discussed for instance in \cite{kushpel1}. See also \cite{kushpel2} and the references therein.

    \subsection{Sobolev Spaces on $\mathbb{M}^d$ and $\mathbb{M}^d\times \mathbb{M}^d$}
    \label{sec:sobolevspaces}
    We now give the definition of Sobolev spaces on $\mathbb{M}^d$ and $\mathbb{M}^d\times \mathbb{M}^d$ which will be essential later to introduce the smoothness properties of the random field $\{X(u), \ u \in \mathbb{M}^d\}$ and its covariance function.
   
   % $$
    %\mathscr{D}_p f = \sum_{\ell \in \Lambda} \sum_{m=1}^{\dime} (1-\lambda_\ell)^{p/2} Y_{\ell,m},
    %$$
    
\begin{definition} \label{def:sobospace}
Let $p \geq 0$. The Sobolev space of order $p$ on $\mathbb{M}^d$ is defined as
\begin{equation*}
    \Hp \!=\!  \Hp (\mathbb{M}^d) \!= \!
    \left\{ f \!\in\! L^2(\mathbb{M}^d) ,\|f\|_{\Hp} < \infty \right\},
  \end{equation*}
  where 
  $$
  \|f\|^2_{\Hp} = \sum_{\ell\in \Lambda} \sum_{m=1}^{\dime} (1-\lambda_\ell)^p  \langle f, Y_{\ell,m} \rangle^2_\rangleL .
  $$
The tensorial Sobolev space of order $p$ on $\mathbb{M}^d \times \mathbb{M}^d$ is
\begin{equation*}
    \Hprod =  \Hprod (\mathbb{M}^d \times \mathbb{M}^d) = \Hp \otimes \Hp \subset L^2 (\mathbb{M}^d \times \mathbb{M}^d).
\end{equation*}
\end{definition}
The parameter $p$ in Definition~\ref{def:sobospace} quantifies the regularity: the higher $p$, the smoother the functions  $f \in \Hp$ and  $g \in \Hprod$. 
Indeed, if we consider the Sobolev operator
    $$
    \mathscr{D} = (\operatorname{Id}-\Delta_{\mathbb{M}^d})^{p/2}, \qquad p \ge 0,
    $$
    we can readily see that
    $$
    \|f\|_\Hp = \|\mathscr{D}f\|_\rangleL, \qquad \|g\|_\Hprod  = \|(\mathscr{D} \otimes \mathscr{D}) g \|_\rangleLL.
    $$          
In the next section, we will only consider Sobolev spaces $\Hp$ and $\Hprod$ that are Reproducing Kernel Hilbert Spaces (RKHS). %~\cite{aronszajn1950theory}
      It turns out that 
    \begin{equation*}
        \Hp \text{ is a RKHS} \quad \Longleftrightarrow \quad  \Hprod \text{ is a RKHS} \quad \Longleftrightarrow  \quad p > d/2   .
    \end{equation*}
    See \cite[Proposition 1]{ejs} and \cite[Corollary 3.4]{gine}.

\section{Covariance Estimation}
In this section, we present the estimation of the covariance function $C$ and the asymptotic theory behind. 

We assume that $\{X(u), \ u \in \mathbb{M}^d\}$ is a random element of $\Hq$, for some $q>d/2$.
For a given $p \ge q$ and $\pen > 0 $, we can define the estimator $C_\pen$ of the covariance function as the solution of the following minimization problem
\begin{equation}\label{eq::cov-est}
\min_{g \in \Hprod}   \sum_{i=1}^n \frac{1}{nr_i(r_i-1)} \sum_{1\le j \ne k \le r_i } (\W_{ij}\W_{ik}  - g(U_{ij}, U_{ik}))^2 + \pen \| g\|^2_\Hprod,
\end{equation}
It is possible to prove that the solution is unique and it is given by
    \begin{equation*}\label{eq:RTformR}
        C_{\eta} =  \sum_{i=1}^n \frac{1}{\sqrt{r_i(r_i-1)}} \sum_{1 \leq j \neq k \leq r_i} \beta^{(\pen)}_{ijk} \,\psi_{\mathscr{D}^* \mathscr{D}} ( \cos \epsilon \rho ( \cdot ,U_{ij}) ) \otimes \psi_{\mathscr{D}^* \mathscr{D}} ( \cos \epsilon \rho ( \cdot ,U_{ik}) )
    \end{equation*}
    for some random $\{\beta^{(\pen)}_{ijk},\ 1 \leq i \leq n, \ 1\leq j \neq k \leq r_i\}$. The function $\psi_{\mathscr{D}^* \mathscr{D}}$ is the \emph{zonal Green's kernel} of $\mathscr{D}^* \mathscr{D}$, defined as
    $$
    \psi_{\mathscr{D}^* \mathscr{D}} ( \cos \epsilon \rho (u,v) ) = \sum_{\ell \in \Lambda} \frac{\kappa_\ell}{(1-\lambda_\ell)^p} P_\ell^{(\alpha,\beta)} (\cos \epsilon \rho (u,v)), \qquad u,v \in \mathbb{M}^d.
    $$
However, the penalty (and hence the solution) can be expressed in terms of any \emph{admissible} operator $\mathscr{D}$ with \emph{spectral growth order} $p$, which induces an equivalent norm -- see \cite[Section 2.2]{ejs}. 
% see also \cite[Definition 3]{ejs}. 

Theorem 1 below establishes a (uniform) rate of convergence for the estimator $C_\eta$, under a suitable condition on the decay of the penalty parameter $\pen$. Such rate is expressed in terms of the number of replicates $n$ and the \emph{average} number of measurement locations $r$, defined as the harmonic mean of the $r_i$'s. The result can thus be interpreted in both the \emph{dense} and \emph{sparse} sampling regimes. In a dense design, $r=r(n)$ is required to diverge with $n$ and it is larger than some order of $n$; on the other hand, in a sparse design, the sampling frequency $r$ is bounded and can be arbitrary small (as small as two). 

In the following, we define the class of probability measures for which our rate is achieved.
\begin{definition}\label{def:pq-cov}
Consider $p>d$, $d/2<q\le p$. Let $\pit$ be the collection of probability measures for $\Hq$-valued processes, such that for any $X$ with probability law $\mathbb{P}_X \in \pit$
$$\mathbb{E}\|X\|^4_{\Hq}\le L, \qquad
\|C\|^2_{\rangleHH} \le K,
$$
for some constants $L, K > 0$. The constants may depend on the manifold $\mathbb{M}^d$.
\end{definition}

\begin{theorem}\label{th::R}
Assume that 
$\mathbb{E}[\epsilon_{11}^4]<\infty$ and the $U_{ij}, i=1,\dots,n, j=1,\dots,r_i$, are independent copies of $U\sim \operatorname{Unif}(\mathbb{M}^d)$. Let $p>d$ and $d/2<q\le p$, and consider the covariance estimator obtained as the minimizer of Equation \eqref{eq::cov-est}. If $\pen \asymp (n r/\log n )^{-2p/(2p+d)}$, then
\begin{align*}
&\lim_{D\to \infty}  \limsup_{n\to\infty} \sup_{\mathbb{P}_X \in \pit}  \mathbb{P}\left (\| C_\pen  - C \|^2_\rangleLL > D\left (\left( \frac{\log n}{nr}\right)^{\frac{2p}{2p+d}} + \frac1n\right ) \right ) \\ &= 0.
\end{align*}
\end{theorem} 

The proof follows the same lines of arguments as \cite[Theorem 4]{ejs}.

%\begin{corollary} \label{coro:asymptoticCetabigO}
%Let $X$ be such that $\mathbb{P}_X \in \pit$. Under the same assumptions of Theorem \ref{th::R},
%\begin{equation*}
  %  \| C_\pen  - C\|^2_\rangleLL = O_\mathbb{P}\left ( \left( \frac{\log n}{n r}\right)^{2p/(2p+d)} + \frac1n\right).
%\end{equation*}
%\end{corollary}

%%%%%%%%%%%%%%%%%%%%%%%%%%%%%%%%%%%%%%%%%%%%%%%%%%%%%%%%%%%%%%%%%%%%%%%%%%%%%%%%

\bigskip
\acknowledgements{The author wishes to thank Maurizia Rossi for insightful discussions on compact two-point homogeneous spaces.}

%%%%%%%%%%%%%%%%%%%%%%%%%%%%%%%%%%%%%%%%%%%%%%%%%%%%%%%%%%%%%%%%%%%%%%%%%%%%%%%%
%%%%%%%%%%%%%%%%%%%%%%%%%%%%%%%%%%%%%%%%%%%%%%%%%%%%%%%%%%%%%%%%%%%%%%%%%%%%%%%%

%%%%%%%%%%%%%%%%%%%%%%%%%%%%%%%%%%%%%%%%%%%%%%%%%%%%%%%%%%%%%%%%%%%%%%%%%%%%%%%%
%%%%%%%%%%%%%%% Can be generated automatically with bibTeX by %%%%%%%%%%%%%%%%%%
%%%%%%%%%%%%%%%%%%%%%%%%%%%%%%%%%%%%%%%%%%%%%%%%%%%%%%%%%%%%%%%%%%%%%%%%%%%%%%%%
%%%%%%%%%%%%%%% \bibliographystyle{abbrv}                     %%%%%%%%%%%%%%%%%%
%%%%%%%%%%%%%%% \bibliography{your-bib-file}                  %%%%%%%%%%%%%%%%%%

\begin{thebibliography}{1}

\bibitem{ejs}
A.~Caponera, J.~Fageot, M.~Simeoni, V.~M.~Panaretos.
\newblock Functional estimation of anisotropic covariance and autocovariance operators on the sphere.
\newblock {\em Electron. J. Stat.}, 16(2):5080--5148, 2022.

\bibitem{gine}
E.~Giné~M..
\newblock The addition formula for the eigenfunctions of the Laplacian.
\newblock {\em Adv. Math.}, 18(1):102--107, 1975.

\bibitem{kushpel1}
A.~Kushpel.
\newblock The Lebesgue constants on projective spaces. 
\newblock {\em Turkish J. Math.}, 45(2):856--863, 2021. 

\bibitem{kushpel2}
A.~Kushpel, S.~A.~Tozoni.
\newblock Entropy and widths of multiplier operators on two-point homogeneous spaces. 
\newblock {\em Constr. Approx.}, 35:137--180, 2012.

\bibitem{wang}
H.-C.~Wang.
\newblock Two-point homogeneous spaces. 
\newblock {\em Ann. Math.}, 55(1):177--191, 1952.


\end{thebibliography}
%%%%%%%%%%%%%%%%%%%%%%%%%%%%%%%%%%%%%%%%%%%%%%%%%%%%%%%%%%%%%%%%%%%%%%%%%%%%%%%%
%%%%%%%%%%%%%%% and then extracted from your-tex-file.bbl     %%%%%%%%%%%%%%%%%%
%%%%%%%%%%%%%%%%%%%%%%%%%%%%%%%%%%%%%%%%%%%%%%%%%%%%%%%%%%%%%%%%%%%%%%%%%%%%%%%%

\end{document}